\documentclass[11pt,a4paper]{article}
\usepackage[T1]{fontenc}
\usepackage[utf8]{inputenc}
\usepackage[english]{babel}
\usepackage{amsmath,amsthm}
\usepackage{newtxtext,newtxmath}
\usepackage[margin=24mm,top=20mm,bottom=22mm,headheight=22pt,headsep=16pt,footskip=24pt]{geometry}
\usepackage{microtype}
\usepackage{graphicx}
\usepackage{caption}
\usepackage{placeins}
\usepackage{xcolor}
\definecolor{ink}{HTML}{242930}
\definecolor{sectionblue}{HTML}{315E77}
\definecolor{questionback}{HTML}{F4F8FC}
\definecolor{questionframe}{HTML}{377EB8}
\definecolor{answerback}{HTML}{F4FAF5}
\definecolor{answerframe}{HTML}{45945B}
\definecolor{aiback}{HTML}{FFF9E9}
\definecolor{aiframe}{HTML}{D99A20}
\definecolor{aititle}{HTML}{F2CD72}
\usepackage{titlesec}
\titleformat{\section}{\large\bfseries\color{sectionblue}}{\thesection}{0.7em}{}
\titleformat{\subsection}{\normalsize\bfseries}{\thesubsection}{0.7em}{}
\titleformat{\subsubsection}{\normalsize\bfseries}{}{0pt}{}
\titlespacing*{\section}{0pt}{2.0ex plus .5ex minus .3ex}{0.8ex}
\titlespacing*{\subsection}{0pt}{1.8ex plus .4ex minus .2ex}{0.6ex}
\titlespacing*{\subsubsection}{0pt}{1.5ex plus .3ex minus .2ex}{0.4ex}
\usepackage{fancyhdr}
\renewcommand{\headrulewidth}{0.3pt}
\renewcommand{\headrule}{\hbox to\headwidth{\color{sectionblue!45}\leaders\hrule height\headrulewidth\hfill}}
\fancypagestyle{plain}{\fancyhf{}\fancyfoot[C]{\thepage}\renewcommand{\headrulewidth}{0pt}}
\usepackage[most]{tcolorbox}
\tcbset{enhanced,breakable,sharp corners,boxrule=0.45pt,
  left=7pt,right=7pt,top=6pt,bottom=6pt,before skip=9pt,after skip=9pt}
\usepackage{enumitem}
\setlist{topsep=4pt,itemsep=2pt,parsep=0pt,leftmargin=1.6em}
\usepackage{needspace}
\usepackage[unicode,colorlinks=true,linkcolor=sectionblue,citecolor=sectionblue,urlcolor=sectionblue]{hyperref}
\hypersetup{pdftitle={Idempotent Ultrafilters and Closed Images of Complete Semilattices},
  pdfsubject={Closed images and natural orders of complete semilattices},
  pdfauthor={Konstantin Kazachenko and Alexander V. Osipov}}
\newtheorem{theorem}{Theorem}[section]
\newtheorem{lemma}[theorem]{Lemma}
\newtheorem{proposition}[theorem]{Proposition}
\theoremstyle{definition}
\newtheorem{question}[theorem]{Question}
\theoremstyle{remark}

\numberwithin{equation}{section}
\DeclareMathOperator{\FP}{FP}
\newcommand{\gen}[1]{\langle #1\rangle}
\newcommand{\cl}[1]{\overline{#1}}
\makeatletter
\renewcommand{\l@section}{\@dottedtocline{1}{0em}{1.4em}}
\renewcommand{\l@subsection}{\@dottedtocline{2}{1.4em}{2.0em}}
\makeatother
\begin{document}
\color{ink}
\thispagestyle{plain}
\begin{center}
{\fontsize{18}{21}\selectfont\bfseries Idempotent Ultrafilters and Closed Images\\[0.15em]
of Complete Semilattices\par}
\vspace{0.8em}
{\normalsize Alexander V. Osipov\footnote{Email: \href{mailto:OAB@list.ru}{\texttt{OAB@list.ru}} .}\quad
Konstantin Kazachenko\footnote{Email: \href{mailto:voice1081@gmail.com}{\texttt{voice1081@gmail.com}}.}\par}
\vspace{0.35em}
{\small\itshape Krasovskii Institute of Mathematics and Mechanics, Yekaterinburg, Russia\par}

\end{center}

\vspace{0.4em}
\begin{list}{}{\setlength{\leftmargin}{14pt}\setlength{\rightmargin}{14pt}}
\item[]\small\textbf{\color{sectionblue}Abstract.}
We prove that every continuous homomorphism from a complete topologized
semilattice to a Hausdorff semitopological semilattice has closed image,
answering a question of Banakh and Bardyla. The main step is the
convergence of idempotent ultrafilters on complete semitopological
semilattices. We also prove that, for every complete subsemilattice
\(X\) of a Hausdorff semitopological semilattice \(Y\), the natural
order of \(X\) is closed in \(Y\times Y\). In particular, the natural
order of every complete Hausdorff semitopological semilattice \(X\)
is closed in \(X\times X\).
\end{list}
\vspace{-0.5em}
\tableofcontents
\section*{Introduction}
\addcontentsline{toc}{section}{Introduction}

The relation between completeness and closedness is a recurring theme in
topological algebra: one seeks intrinsic conditions ensuring that a
structure, or each of its continuous homomorphic images, is closed in
any suitable Hausdorff extension. For semilattices, this question brings
together the topology and the natural order.

Stepp characterized absolutely maximal discrete semilattices
\cite[Theorem~9]{S75}. In modern terminology, a discrete semilattice
has this property if and only if it is chain-finite, meaning that every
chain is finite. In particular, every homomorphism from a chain-finite
semilattice into a Hausdorff topological semilattice has closed image.
Gutik and Repov\v{s} characterized \(H\)-closed linearly ordered
topological semilattices by order completeness and closure conditions
at suprema and infima, and proved that every such semilattice is
absolutely \(H\)-closed~\cite[Theorems~2 and~3]{GR08}.

Banakh and Bardyla developed a general theory of chain-compact and
chain-finite semilattices~\cite{BB19}. They proved that every continuous
homomorphism from a complete topologized semilattice into a Hausdorff
topological semilattice has closed image. They also proved that every
homomorphism from a chain-finite semilattice into a Hausdorff
semitopological semilattice has closed image. These results led to the
question whether the first conclusion remains valid when the target
operation is only separately continuous. The question is discussed in
Banakh and Bardyla's survey \emph{Complete topologized posets and
semilattices}~\cite{BB}.

\begin{figure}[!htbp]
\centering
\includegraphics[width=\linewidth]{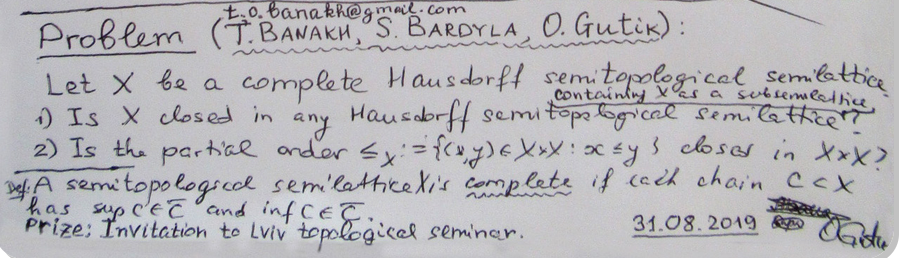}
\caption{The entry by T.~Banakh, S.~Bardyla, and O.~Gutik in the
\emph{Lviv Scottish Book}, volume~2, p.~161, dated 31~August~2019.}
\label{fig:scottish-book}
\end{figure}

The closedness question, together with the question of whether the
natural order is closed in \(X\times X\), was also recorded by Taras
Banakh, Serhii Bardyla, and Oleg Gutik in the modern \emph{Lviv Scottish
Book}, volume~2, p.~161, in an entry dated 31~August~2019~\cite{BBG19}.
Figure~\ref{fig:scottish-book} reproduces this handwritten entry.

Several positive cases were established. Banakh, Bardyla, and Ravsky
proved the closed-image conclusion for functionally Hausdorff
semitopological targets, and more generally under a weaker separation
axiom~\cite{BBR19}. Banakh and Bardyla proved it for sequential
Hausdorff semitopological targets~\cite{BB20}. Banakh, Bardyla, and
Gutik obtained it for \(\omega\)-Lawson semitopological
semilattices~\cite{BBG21}.

We prove the closed-image assertion for arbitrary Hausdorff
semitopological targets. The main step is to show that every idempotent
ultrafilter on a complete semitopological semilattice converges.
We use the standard semigroup structure on the Stone--\v{C}ech
compactification of a discrete semigroup~\cite{HS}.

We also show that the natural order of a complete Hausdorff
semitopological semilattice \(X\), regarded as a subset of
\(X\times X\), is closed. More generally, if \(X\) is a complete
subsemilattice of a Hausdorff semitopological semilattice \(Y\), then
\(\{(x,y)\in X\times X:xy=x\}\) is closed in \(Y\times Y\).
This latter formulation, considered under additional assumptions
in~\cite{BBR19,BBG21}, is equivalent to the closed-image assertion
when both are required for all such semilattices.

Section~\ref{sec:problem} proves convergence of idempotent ultrafilters
and the closed-image theorem. Section~\ref{sec:natural-order}
establishes closedness of the natural order and the equivalence just
mentioned.

\FloatBarrier
\section{Closed images of complete semilattices}\label{sec:problem}
For a semigroup \(T\) and each \(t\in T\), let
\[
l_t:T\to T,\quad l_t(x)=tx,\qquad
r_t:T\to T,\quad r_t(x)=xt
\]
denote the left and right shifts. A \emph{semilattice} is a commutative
semigroup consisting of idempotents. Its natural partial order is given by
\[
x\le y\quad\Longleftrightarrow\quad xy=x.
\]
A \emph{topologized semilattice} is a semilattice endowed with a topology.
It is \emph{semitopological} if every shift is continuous. Following
Banakh and Bardyla~\cite{BB}, it is \emph{complete} if each nonempty
chain \(C\) has an infimum and a supremum in the semilattice and
\[
\inf C\in\cl C,\qquad \sup C\in\cl C.
\]

\begin{tcolorbox}[colback=questionback,colframe=questionframe]
\begin{question}[Banakh and Bardyla~\cite{BB}]\label{q:main}
Let \(h:X\to Y\) be a continuous homomorphism from a complete topologized
semilattice \(X\) to a Hausdorff semitopological semilattice \(Y\).
Must the image \(h(X)\) be closed in \(Y\)?
\end{question}
\end{tcolorbox}

The answer is affirmative, as follows from
Theorem~\ref{thm:closed-subsemilattice} below.

\subsection{Idempotent ultrafilters}\label{sec:solution}
Except where stated otherwise, semilattices are nonempty, and closures
are taken in the specified topology. For a subset \(A\) of a
semilattice \(S\), write
\[
\gen A=\{a_1\cdots a_n:n\ge1,\ a_1,\ldots,a_n\in A\}.
\]
We put \(\gen\varnothing=\varnothing\). For every infinite cardinal
\(\kappa\), if \(|A|\le\kappa\), then \(|\gen A|\le\kappa\).

We shall use the fact that the closure of a subsemilattice in a
semitopological semilattice is a subsemilattice
\cite[Theorem~2.27]{HS}. We shall also use two completeness properties
appearing in the closed-image proof in~\cite{BB19}: every nonempty
closed subsemilattice of a complete topologized semilattice has a least
element, and every decreasing ordinal-indexed family of nonempty
closed subsemilattices has nonempty intersection.

Let \(S_d\) denote the semilattice \(S\) with the discrete topology.
Identify \(\beta S_d\) with the compact Hausdorff space of ultrafilters
on \(S\). For each \(A\subseteq S\), put
\[
\widehat A=\{p\in\beta S_d:A\in p\}.
\]
These sets form a clopen base of \(\beta S_d\).
The standard extension of multiplication to \(\beta S_d\) makes it a
compact Hausdorff right-topological semigroup
\cite[Theorems~4.1 and~4.4]{HS}. For every \(p,q\in\beta S_d\) and
every \(A\subseteq S\), the product is characterized by
\begin{equation}\label{eq:ultrafilter-product}
A\in pq\quad\Longleftrightarrow\quad
\{s\in S:l_s^{-1}(A)\in q\}\in p
\end{equation}
\cite[Theorem~4.12(b)]{HS}.

For an idempotent \(p\in\beta S_d\) and a subset \(A\subseteq S\), set
\[
A^\star=A\cap\{s\in S:l_s^{-1}(A)\in p\}.
\]
The notation \(A^\star\) refers to the fixed ultrafilter \(p\).
By \cite[Definition~4.13 and Lemma~4.14]{HS}, if \(A\in p\), then
\(A^\star\in p\), and for every \(s\in A^\star\),
\begin{equation}\label{eq:star-stability}
l_s^{-1}(A^\star)\in p.
\end{equation}
For a sequence \((x_n)_{n<\omega}\) and every \(n<\omega\), write
\[
\FP(x_j:j\ge n)=
\left\{\prod_{j\in F}x_j:
\varnothing\ne F\subseteq\{n,n+1,\ldots\},\ F\text{ finite}\right\}.
\]

\begin{lemma}\label{lem:diagonal}
Let \(p\in\beta S_d\) be idempotent, and let \((V_n)_{n<\omega}\) be a
decreasing sequence of members of \(p\). There exists a sequence
\((x_n)_{n<\omega}\) in \(S\) such that, for every nonempty finite
set \(F\subseteq\omega\),
\begin{equation}\label{eq:finite-invariant}
\prod_{j\in F}x_j\in V_{\min F}^\star.
\end{equation}
\end{lemma}

\begin{proof}
We construct \((x_n)\) recursively so that \eqref{eq:finite-invariant}
holds for every nonempty subset of the indices already chosen.
Choose \(x_0\in V_0^\star\).

Suppose that \(x_0,\ldots,x_{n-1}\) have been chosen, where \(n\ge1\).
Choose
\[
x_n\in V_n^\star\cap
\bigcap_{\varnothing\ne F\subseteq\{0,\ldots,n-1\}}
l_{\prod_{j\in F}x_j}^{-1}\bigl(V_{\min F}^\star\bigr).
\]
Every set in this finite intersection belongs to \(p\), by
\eqref{eq:star-stability} and the induction hypothesis. The intersection
is therefore nonempty. The choice gives \(x_n\in V_n^\star\) and,
for every nonempty \(F\subseteq\{0,\ldots,n-1\}\),
\[
\left(\prod_{j\in F}x_j\right)x_n\in V_{\min F}^\star.
\]
Together with the induction hypothesis, these inclusions establish
\eqref{eq:finite-invariant} for all nonempty subsets of
\(\{0,\ldots,n\}\). Every finite subset of \(\omega\) is covered
at some stage, so the resulting sequence has the required property.
\end{proof}

\begin{proposition}\label{prop:simultaneous}
Let \(S\) be a complete semitopological semilattice, let
\(p\in\beta S_d\) be idempotent, and let \((C_i)_{i\in I}\) be a family
of closed subsets of \(S\) belonging to \(p\).
For every infinite cardinal \(\kappa\) and every \(R\subseteq S\) with
\(|I|,|R|\le\kappa\), there exist a subsemilattice \(D\subseteq S\) and
a point \(b\in S\) such that
\begin{equation}\label{eq:support}
R\subseteq D,\qquad |D|\le\kappa,
\end{equation}
and, for every \(i\in I\),
\begin{equation}\label{eq:pure}
b\in\cl{C_i^\star\cap D}.
\end{equation}
For every \(i\in I\) and every \(r\in R\cap C_i^\star\), we also have
\begin{equation}\label{eq:weighted}
rb\in\cl{C_i^\star\cap D}.
\end{equation}
\end{proposition}

\begin{proof}
We argue by induction on infinite cardinals. First suppose that
\(\kappa=\aleph_0\). The family
\[
\{S\}\ \cup\ \{C_i^\star:i\in I\}\ \cup\
\{l_r^{-1}(C_i^\star):i\in I,\ r\in R\cap C_i^\star\}
\]
is nonempty and at most countable, and all its members belong to \(p\)
by \eqref{eq:star-stability}. Enumerate it as \((W_n)_{n<\omega}\),
with repetitions if necessary, and put \(V_n=\bigcap_{j\le n}W_j\).
Choose \((x_n)_{n<\omega}\) by Lemma~\ref{lem:diagonal}. For every
\(m<\omega\) and every nonempty finite
\(F\subseteq\{m,m+1,\ldots\}\), the lemma gives
\[
\prod_{j\in F}x_j\in V_{\min F}^\star
\subseteq V_{\min F}\subseteq V_m.
\]
Thus \(\FP(x_j:j\ge m)\subseteq V_m\) for every \(m<\omega\).

Let \(D=\gen{R\cup\{x_n:n<\omega\}}\), which satisfies
\eqref{eq:support}. For every \(m<\omega\) and every \(n\ge m\), put
\[
u_{m,n}=x_mx_{m+1}\cdots x_n\in V_m\cap D,
\]
and let \(c_m=\inf_{n\ge m}u_{m,n}\). For fixed \(m\), the elements
\(u_{m,n}\) form a decreasing chain, so completeness gives
\begin{equation}\label{eq:row-closure}
c_m\in\cl{\{u_{m,n}:n\ge m\}}.
\end{equation}
For every \(n\ge m+1\), we have
\(c_m\le u_{m,n}\le u_{m+1,n}\). Hence \(c_m\le c_{m+1}\).
By completeness, \(b=\sup_{m<\omega}c_m\) exists and belongs to the
closure of every tail of \((c_m)\).

Fix \(i\in I\), and choose \(j<\omega\) such that \(W_j=C_i^\star\).
For every \(m\ge j\) and every \(n\ge m\),
\[
u_{m,n}\in V_m\cap D\subseteq C_i^\star\cap D.
\]
By \eqref{eq:row-closure}, \(c_m\in\cl{C_i^\star\cap D}\) for every
\(m\ge j\), and therefore
\[
b\in\cl{\{c_m:m\ge j\}}\subseteq\cl{C_i^\star\cap D}.
\]

Fix \(r\in R\cap C_i^\star\), and choose \(k<\omega\) such that
\(W_k=l_r^{-1}(C_i^\star)\). For every \(m\ge k\) and every \(n\ge m\),
\[
u_{m,n}\in V_m\subseteq l_r^{-1}(C_i^\star),
\qquad ru_{m,n}\in C_i^\star\cap D.
\]
Continuity of \(l_r\) and \eqref{eq:row-closure} imply that, for every
\(m\ge k\),
\[
rc_m\in\cl{\{ru_{m,n}:n\ge m\}}\subseteq\cl{C_i^\star\cap D}.
\]
Applying \(l_r\) to the closure of the tail \(\{c_m:m\ge k\}\), we get
\[
rb\in\cl{\{rc_m:m\ge k\}}\subseteq\cl{C_i^\star\cap D}.
\]
This proves the countable case.

Let \(\kappa>\aleph_0\), and suppose that the proposition holds for
every smaller infinite cardinal. Fix a complete semitopological
semilattice \(S\), an idempotent \(p\in\beta S_d\), a family
\((C_i)_{i\in I}\) of closed members of \(p\), and a subset
\(R\subseteq S\) with \(|I|,|R|\le\kappa\). Choose injections
\(\iota:I\to\kappa\) and \(\rho:R\to\kappa\). For every \(\alpha<\kappa\), put
\[
I_\alpha=\{i\in I:\iota(i)\le\alpha\},\qquad
R_\alpha^0=\{r\in R:\rho(r)\le\alpha\},\qquad
\lambda_\alpha=\max\{\aleph_0,|\alpha|\}.
\]
We construct subsemilattices \((D_\alpha)_{\alpha<\kappa}\) of \(S\)
and points \((b_\alpha)_{\alpha<\kappa}\) of \(S\) by transfinite
recursion. At each stage, write
\[
R_\alpha=\gen{\bigcup_{\beta<\alpha}D_\beta\ \cup\ R_\alpha^0}.
\]
We require
\begin{equation}\label{eq:stage-support}
R_\alpha\subseteq D_\alpha,\qquad |D_\alpha|\le\lambda_\alpha,
\end{equation}
and, for every \(i\in I_\alpha\),
\begin{equation}\label{eq:stage-pure}
b_\alpha\in\cl{C_i^\star\cap D_\alpha}.
\end{equation}
For every \(i\in I_\alpha\) and every \(q\in R_\alpha\cap C_i^\star\),
we also require
\begin{equation}\label{eq:stage-weighted}
q b_\alpha\in\cl{C_i^\star\cap D_\alpha}.
\end{equation}

At stage \(0\), we have \(R_0=\gen{R_0^0}\) and
\(\lambda_0=\aleph_0\). The sets \(I_0\) and \(R_0^0\) each have
at most one element. The countable case, applied to
\((C_i)_{i\in I_0}\) and \(R_0\), therefore supplies \(D_0,b_0\)
satisfying \eqref{eq:stage-support}--\eqref{eq:stage-weighted}.

Fix \(0<\alpha<\kappa\), and suppose that \(D_\beta,b_\beta\) have
been constructed for every \(\beta<\alpha\). Since
\(\lambda_\beta\le\lambda_\alpha\) for every \(\beta<\alpha\),
\[
\left|\bigcup_{\beta<\alpha}D_\beta\right|
\le |\alpha|\cdot\lambda_\alpha\le\lambda_\alpha.
\]
Also \(|I_\alpha|,|R_\alpha^0|\le|\alpha+1|\le\lambda_\alpha\),
so \(|R_\alpha|\le\lambda_\alpha<\kappa\). Apply the induction
hypothesis at \(\lambda_\alpha\) to the family
\((C_i)_{i\in I_\alpha}\) and the set \(R_\alpha\). It gives
\(D_\alpha,b_\alpha\) satisfying
\eqref{eq:stage-support}--\eqref{eq:stage-weighted}, completing the recursion.

For every \(\beta<\alpha<\kappa\),
\(D_\beta\subseteq R_\alpha\subseteq D_\alpha\). Consequently,
\[
D=\bigcup_{\alpha<\kappa}D_\alpha
\]
is a subsemilattice of cardinality at most \(\kappa\).
For every \(r\in R\),
\(r\in R_{\rho(r)}^0\subseteq D_{\rho(r)}\subseteq D\), so
\eqref{eq:support} holds.

Fix \(\beta<\alpha<\kappa\) and \(i\in I_\alpha\).
Since \(D_\beta\subseteq R_\alpha\), condition
\eqref{eq:stage-weighted} gives
\[
r_{b_\alpha}(C_i^\star\cap D_\beta)
\subseteq\cl{C_i^\star\cap D_\alpha}.
\]
Continuity of \(r_{b_\alpha}\) yields
\begin{equation}\label{eq:bootstrap}
r_{b_\alpha}\bigl(\cl{C_i^\star\cap D_\beta}\bigr)
\subseteq\cl{C_i^\star\cap D_\alpha}.
\end{equation}

For every \(i\in I\) and every nonempty finite increasing sequence
\(\alpha_1<\cdots<\alpha_n<\kappa\) with \(\iota(i)\le\alpha_1\), we claim that
\begin{equation}\label{eq:local-pure}
b_{\alpha_1}\cdots b_{\alpha_n}\in\cl{C_i^\star\cap D_{\alpha_n}}.
\end{equation}
We prove this by induction on \(n\). For \(n=1\), it is
\eqref{eq:stage-pure}. For \(n>1\), the product of the first \(n-1\)
factors belongs to \(\cl{C_i^\star\cap D_{\alpha_{n-1}}}\), and
\eqref{eq:bootstrap} with \(\beta=\alpha_{n-1}\), \(\alpha=\alpha_n\)
gives the required inclusion.

For every \(i\in I\), every \(r\in R\cap C_i^\star\), and every
nonempty finite increasing sequence \(\alpha_1<\cdots<\alpha_n<\kappa\)
with \(\max\{\iota(i),\rho(r)\}\le\alpha_1\), we similarly obtain
\begin{equation}\label{eq:local-weighted}
r b_{\alpha_1}\cdots b_{\alpha_n}\in\cl{C_i^\star\cap D_{\alpha_n}}.
\end{equation}
For \(n=1\), this follows from \eqref{eq:stage-weighted}, because
\(r\in R_{\alpha_1}^0\subseteq R_{\alpha_1}\).
The induction step again follows from \eqref{eq:bootstrap}.

For each \(\gamma<\kappa\), put
\[
P_\gamma=\gen{\{b_\alpha:\gamma\le\alpha<\kappa\}},\qquad
E_\gamma=\cl{P_\gamma}.
\]
The sets \(E_\gamma\) form a decreasing family of nonempty closed
subsemilattices, by \cite[Theorem~2.27]{HS}. The intersection property
recalled above gives a point \(b\in\bigcap_{\gamma<\kappa}E_\gamma\).

Fix \(i\in I\) and \(\gamma\ge\iota(i)\). Idempotence allows
repeated factors in a finite product to be removed, so
\eqref{eq:local-pure} gives
\(P_\gamma\subseteq\cl{C_i^\star\cap D}\). Hence
\(E_\gamma\subseteq\cl{C_i^\star\cap D}\), and \eqref{eq:pure} follows.

Fix \(r\in R\cap C_i^\star\) and
\(\gamma\ge\max\{\iota(i),\rho(r)\}\). By \eqref{eq:local-weighted},
\(rP_\gamma\subseteq\cl{C_i^\star\cap D}\). Continuity of \(l_r\)
therefore gives
\[
rE_\gamma\subseteq\cl{rP_\gamma}\subseteq\cl{C_i^\star\cap D}.
\]
Since \(b\in E_\gamma\), this proves \eqref{eq:weighted} and completes
the induction.
\end{proof}

\begin{theorem}\label{thm:idempotent}
Every idempotent ultrafilter on a complete semitopological semilattice
converges.
\end{theorem}

\begin{proof}
Let \(S\) be a complete semitopological semilattice and let
\(p\in\beta S_d\) be idempotent. Apply
Proposition~\ref{prop:simultaneous} with \(R=\varnothing\) to the
family of all closed members of \(p\), choosing an infinite cardinal
at least as large as this family. We obtain a point \(b\in S\) such that
\[
b\in\cl{C^\star\cap D}\subseteq C
\]
for every closed \(C\in p\). If an open neighbourhood \(U\) of
\(b\) did not belong to \(p\), then the closed set \(S\setminus U\)
would belong to \(p\), a contradiction. Thus \(p\to b\).
\end{proof}

\subsection{Closedness of images}\label{sec:images}
We use the Ellis--Numakura lemma: every nonempty compact Hausdorff
right-topological semigroup contains an idempotent
\cite[Theorem~2.5]{HS}.

\begin{lemma}\label{lem:boundary}
Let \(S\) be a nonempty subsemilattice of a semitopological semilattice
\(Y\), and let \(y\in\cl S^{\,Y}\). Then
\[
K_y=\{p\in\beta S_d:p\to y\text{ in }Y\}
\]
is a nonempty compact subsemigroup of \(\beta S_d\).
\end{lemma}

\begin{proof}
Let \(\mathcal N_y\) be the family of open neighbourhoods of \(y\) in \(Y\).
Then
\[
K_y=\bigcap_{U\in\mathcal N_y}\widehat{U\cap S}.
\]
The sets \(U\cap S\), with \(U\in\mathcal N_y\), have the finite
intersection property because \(y\in\cl S^{\,Y}\). Hence the clopen
sets \(\widehat{U\cap S}\) have the finite intersection property in
\(\beta S_d\), and their intersection is nonempty and compact.

Fix \(p,q\in K_y\) and \(U\in\mathcal N_y\).
Since \(y^2=y\) and \(r_y:Y\to Y\) is continuous,
\[
B=\{s\in S:sy\in U\}\in p.
\]
For every \(s\in B\), the set \(\{t\in Y:st\in U\}\) is an open
neighbourhood of \(y\). Since \(q\to y\), its trace on \(S\) belongs
to \(q\), so \(l_s^{-1}(U\cap S)\in q\), where \(l_s:S\to S\).
Thus
\[
B\subseteq\{s\in S:l_s^{-1}(U\cap S)\in q\}.
\]
By \eqref{eq:ultrafilter-product}, \(U\cap S\in pq\).
As \(U\) was arbitrary, \(pq\in K_y\).
\end{proof}

By \cite[Lemma~5.3]{BB}, a continuous homomorphic image of a complete
topologized semilattice in a Hausdorff semitopological semilattice is
complete. Consequently, the closed-image assertion in
Question~\ref{q:main} is equivalent to the following theorem: one
direction follows by applying the theorem to \(h(X)\), and the other
by taking inclusion homomorphisms.

\begin{theorem}\label{thm:closed-subsemilattice}
Every complete subsemilattice of a Hausdorff semitopological semilattice
is closed.
\end{theorem}

\begin{proof}
Let \(S\) be a nonempty complete subsemilattice of a Hausdorff
semitopological semilattice \(Y\), and fix \(y\in\cl S^{\,Y}\).
By Lemma~\ref{lem:boundary} and the Ellis--Numakura lemma, \(K_y\)
contains an idempotent \(p\). By Theorem~\ref{thm:idempotent}, \(p\)
converges to a point \(b\in S\). The inclusion \(S\to Y\) is
continuous, so \(p\) converges to both \(b\) and \(y\) in \(Y\).
Hausdorffness gives \(b=y\), and therefore \(y\in S\).
\end{proof}

\section{Closedness of the natural order}\label{sec:natural-order}
Products of complete topologized semilattices are complete
\cite[Proposition~5.2]{BB}. We apply this fact to the natural order
regarded as a subsemilattice of the square.

\begin{proposition}\label{prop:complete-order-graph}
Let \(X\) be a complete topologized semilattice. Then
\[
P_X=\{(x,y)\in X\times X:x\le y\}
\]
is a complete topologized subsemilattice of \(X\times X\).
\end{proposition}

\begin{proof}
If \((x,y),(u,v)\in P_X\), then
\((xu)(yv)=(xy)(uv)=xu\), so \((xu,yv)\in P_X\).
Let \(C\subseteq P_X\) be a nonempty chain. Completeness of \(X^2\)
gives its infimum and supremum
\[
(a,b)=\bigl(\inf\pi_1(C),\inf\pi_2(C)\bigr),\qquad
(c,d)=\bigl(\sup\pi_1(C),\sup\pi_2(C)\bigr),
\]
both belonging to \(\cl C^{\,X^2}\).
For every \((x,y)\in C\), we have \(a\le x\le y\), so
\(a\le b\). Likewise, \(x\le y\le d\) for every \((x,y)\in C\),
so \(c\le d\). Thus both pairs belong to \(P_X\), and they are
also the infimum and supremum of \(C\) in \(P_X\). Since
\(\cl C^{\,P_X}=P_X\cap\cl C^{\,X^2}\), the semilattice \(P_X\)
is complete.
\end{proof}

\begin{theorem}\label{thm:closed-natural-order}
Let \(X\) be a complete subsemilattice of a Hausdorff semitopological
semilattice \(Y\). Then
\[
P_X=\{(x,y)\in X\times X:xy=x\}
\]
is closed in \(Y\times Y\). In particular, the natural order of every
complete Hausdorff semitopological semilattice \(X\) is closed in
\(X\times X\).
\end{theorem}

\begin{proof}
By Proposition~\ref{prop:complete-order-graph}, \(P_X\) is a complete
subsemilattice of \(X^2\), and hence of \(Y^2\), with the same subspace
topology. The product \(Y^2\) is a Hausdorff semitopological semilattice.
Theorem~\ref{thm:closed-subsemilattice}, applied to \(P_X\subseteq Y^2\),
therefore shows that \(P_X\) is closed in \(Y^2\).
\end{proof}

Conversely, closedness of \(P_X\) in \(Y^2\) implies closedness of
\(X\) in \(Y\): for the continuous diagonal map
\(\delta:Y\to Y^2\), \(\delta(y)=(y,y)\), we have
\(\delta^{-1}(P_X)=X\). Thus, Theorems~\ref{thm:closed-natural-order} and
\ref{thm:closed-subsemilattice} are equivalent in the class of Hausdorff
semitopological semilattices. By the
image-preservation result cited above, they are also equivalent to the
closed-image assertion in Question~\ref{q:main}.

\section*{Acknowledgments}
The research was supported by the Russian Science Foundation (RSF Grant No.~26-21-00196).

\section*{AI best practices}
\addcontentsline{toc}{section}{AI best practices}
Codex assisted with the development and revision of the argument.
Parallel agents were used to examine intermediate proofs.

\end{document}